\documentclass[11pt,leqno]{amsart}
\usepackage{amssymb,verbatim,enumerate,ifthen,times}
\usepackage[mathscr]{eucal}
\usepackage[utf8]{inputenc}
\usepackage[T1]{fontenc}
\usepackage{dsfont}
\usepackage{accents}

\makeatletter
\@namedef{subjclassname@2020}{\textup{2020} Mathematics Subject Classification}
\makeatother

\usepackage{accents}
\def\N{\mathbb{N}}

\def\R{\mathbb{R}}

\def\diam{\mbox{\rm diam}}
\def\Ref(#1|#2){(#1\hspace{.3mm}|\hspace{.3mm}#2)}
\newtheorem{theorem}{Theorem}
\newtheorem*{theorem*}{Theorem}
\long\def\Thm#1#2{\ifthenelse{\equal{#1}{*}}{\begin{theorem*}#2\end{theorem*}}
             {\begin{theorem}\label{T#1}#2\end{theorem}}}
\newtheorem{Atheorem}{Theorem}

\newtheorem{proposition}[theorem]{Proposition}
\newtheorem*{proposition*}{Proposition}
\long\def\Prp#1#2{\ifthenelse{\equal{#1}{*}}{\begin{proposition*}#2\end{proposition*}}
             {\begin{proposition}\label{P#1}#2\end{proposition}}}
\def\prp#1{Proposition~\ref{P#1}}

\newtheorem{corollary}[theorem]{Corollary}
\newtheorem*{corollary*}{Corollary}
\long\def\Cor#1#2{\ifthenelse{\equal{#1}{*}}{\begin{corollary*}#2\end{corollary*}}
             {\begin{corollary}\label{C#1}#2\end{corollary}}}

\newtheorem{lemma}[theorem]{Lemma}
\newtheorem*{lemma*}{Lemma}
\long\def\Lem#1#2{\ifthenelse{\equal{#1}{*}}{\begin{lemma*}#2\end{lemma*}}
             {\begin{lemma}\label{L#1}#2\end{lemma}}}
\def\lem#1{Lemma~\ref{L#1}}

\theoremstyle{definition}
\newtheorem{definition}[theorem]{Definition}
\newtheorem*{definition*}{Definition}
\long\def\Defn#1#2{\ifthenelse{\equal{#1}{*}}{\begin{definition*}\rm #2\end{definition*}}
             {\begin{definition}\label{D#1}\rm #2\end{definition}}}

\newtheorem{remark}[theorem]{Remark}
\newtheorem*{remark*}{Remark}
\long\def\Rem#1#2{\ifthenelse{\equal{#1}{*}}{\begin{remark*}\rm #2\end{remark*}}
             {\begin{remark}\label{R#1}\rm #2\end{remark}}}

\newtheorem{example}{Example}
\newtheorem*{example*}{Example}
\long\def\Exa#1#2{\ifthenelse{\equal{#1}{*}}{\begin{example*}\rm #2\end{example*}}
             {\begin{example}\label{Ex#1}\rm #2\end{example}}}

\def\eq#1{{\rm(\ref{E#1})}}
\def\Eq#1#2{\ifthenelse{\equal{#1}{*}}
  {\begin{equation*}\begin{aligned}[]#2\end{aligned}\end{equation*}}
  {\begin{equation}\begin{aligned}\label{E#1}#2\end{aligned}\end{equation}}}

\begin{document}

\date{\today}

\title{On the $\sigma$-balancing property of multivariate generalized quasi-arithmetic means}  
\author[T. Kiss]{Tibor Kiss}
\author[G. Nagy]{Gergő Nagy}
\email{kiss.tibor@science.unideb.hu}
\email{nagyg@science.unideb.hu}
\address{Institute of Mathematics, University of Debrecen, P. O. Box 400, H-4002 Debrecen, Hungary and HUN-REN-UD Equations, Functions, Curves 
and their Applications Research Group}
\thanks{This project has received funding from the HUN-REN Hungarian Research Network. The research of the first author was supported in part by NKFIH 
Grant K-134191.}

\begin{abstract}
The aim of this paper is to characterize the so-called $\sigma$-balancing property in the class of generalized quasi-arithmetic means. In general, the 
question is 
whether those elements of a given family of means that possess this property are quasi-arithmetic.

The first result in the latter direction is due to G. Aumann who showed that a balanced complex mean is necessariliy quasi-arithmetic provided that it 
is analytic. Then Aumann characterized quasi-arithmetic means among Cauchy means in terms of the balancing property. These results date back to the 
1930s. In 2015, Lucio R. Berrone, generalizing balancedness, concluded that a mean having that more general property is quasi-arithmetic if it is 
symmetric, strict and continuously differentiable. A common feature of these results is that they assume a certain order of differentiability of the 
mean whether or not it is a natural condition.

In 2020, the balancing property was characterized in the family of generalized quasi-arithmetic means of two variables under only natural conditions, 
namely continuity and strict monotonicity of their generating functions. Here we extend the corresponding result for multivariate generalized 
quasi-arithmetic means by relaxing the conditions on the generating functions and considering the more general $\sigma$-balancing property.

\end{abstract}

\maketitle

\section{Introduction}

We introduce some concepts and notation that will be used in the paper. We fix a number $n\in\N$ with $n\geq 2$ and a non-trivial interval 
$I\subseteq\R$, i.e., one of positive length. The interior of $I$ will be denoted by $I^\circ$. In general, a real-valued function 
$M:I^n\to\R$ is called a \emph{mean of $n$ variables on $I$} if for all vectors $(x_1,\dots,x_n)\in I^n$, the chain of inequalities
\Eq{mean}{
	\min(x_1,\dots,x_n)\leq M(x_1,\dots,x_n)\leq\max(x_1,\dots,x_n)
}
holds true. $M$ is called \emph{strict} if the relations in 
\eq{mean} are strict whenever $x_1,\dots,x_n$ are pairwise different.  We also note that each mean is a \emph{reflexive} function, that is, 
$M(t,\dots,t)=t$ holds for all $t\in I$. Finally, $M$ is called 
\emph{symmetric} if for any given bijection $\sigma:\{1,\dots,n\}\to\{1,\dots,n\}$, we have that $M(x_1,\dots,x_n)=M(x_{\sigma(1)},\dots,x_{\sigma(n)})$.

In order to define the property investigated in the paper precisely, we introduce the following notation. For $t\in I$, $x=(x_1,\dots,x_n)\in I^n$ 
and $k\in\{1,\dots,n\}$, define the vector $u_k(x,t)\in I^n$ by
\Eq{*}{
	u_k(x,t)_i:=
	\begin{cases}
		t,&\text{if }i=k\\
		x_i,&\text{if }i\neq k,
	\end{cases}
}
where $i=1,\ldots,n$. Then, for a bijection $\sigma:\{1,\dots,n\}\to\{1,\dots,n\}$, we say that a mean $M:I^n\to\R$ of $n$ variables is 
\emph{$\sigma$-balanced} if for all $x=(x_1,\dots,x_n)\in I^n$, the value $u=M(x)$ fulfills the 
fixed point equation
\Eq{Bal}{
M\big(M(u_{\sigma(1)}(x,u)),\dots,M(u_{\sigma(n)}(x,u))\big)=u.
}
On the one hand, originally, the balancing property of means was defined only in the case where $\sigma(k)=n-k+1$ for $k=1,\ldots,n$. The 
$\sigma$-balanced means corresponding to this permutation will be called balanced means. On the other hand, we think that the most natural notion of 
that property is the one in which $\sigma$ is the identity. The previous considerations motivated the definition of $\sigma$-balanced means.
 
In the case where $n=2$, $\sigma(1)=2$, and $\sigma(2)=1$, the relation \eq{Bal} reduces to \emph{Aumann's Equation}
\Eq{*}{
	M(M(s,M(s,t)),M(M(s,t),t))=M(s,t)\qquad(s,t\in I),
}
which was investigated in \cite{Aum35}.

It is easy to see that, for any bijection $\sigma$, the arithmetic mean is $\sigma$-balanced. Indeed, since it is symmetric, during the 
proof of this claim, we may and do assume that $\sigma$ is the identity. Now, if $M(x_1,\dots,x_n)=\frac1n(x_1+\dots+x_n)=:u$, then 
\Eq{*}{
	M\big(M(u_1(x,u)),\dots,M(u_n(x,u))\big)
	&=\frac1n\sum\limits_{k=1}^nM(u_k(x,u))=\frac1{n^2}\sum\limits_{k=1}^n\Big(u+\sum\limits_{j=1,\ j\neq k}^nx_j\Big)\\
	&=\frac1{n^2}\sum\limits_{k=1}^n(u+nu-x_k)=\frac1{n^2}(nu+n^2u-nu)=u.
}

A much more general example is the class of quasi-arithmetic means. These give a continuous, symmetric, strict solution to \eq{Bal}. We know of no 
other solution of this type. In fact, the problem of the existence of a continuous, symmetric, strict solution of \eq{Bal}, which fails to be a 
quasi-arithmetic mean is an open question. It is also worth noting that the balancing property is independent of these three properties in some sense.

In general, it is an interesting question whether having a mean, beside the $\sigma$-balancing property, what we need to assume to conclude that it is 
quasi-arithmetic. The first remarkable investigations in this direction are due to Georg Aumann \cite{Aum35,Aum37}. In 1935, considering the problem 
on the complex plane, in the paper \cite{Aum35}, it was shown that among analytic means (i.e. reflexive, symmetric, holomorphic functions defined on 
an open ball in $\mathbb{C}^n$), the only balanced ones are the analytic quasi-arithmetic means. Then, two years later, Aumann proved that the 
balancing property characterizes quasi-arithmetic means among Cauchy means \cite{Aum37}. The reader will find several balanced means with different 
regularity properties in \cite{Kis20} and \cite{Ber15}.

In 2015, Lucio R. Berrone \cite{Ber15} introduced and investigated a generalized version of the balancing property. Under various combinations of the 
conditions symmetry, strictness, and continuous differentiability, he obtained that the means having this variant of the latter property are 
quasi-arithmetic ones.

A common feature of these results is that they assume some degree of a differentiability property about the means or their generating functions. In 
2020, under only natural regularity conditions, the first author characterized balanced generalized quasi-arithmetic means of two variables 
\cite{Kis20}. The corresponding theorem is analogous to Aumann's result on Cauchy means, namely the solution is the class of quasi-arithmetic means of 
two variables. In this paper, we want to generalize this statement in several directions. We extend the result to the $n$ variable case, we require 
more relaxed conditions on the generating functions and we consider the more general $\sigma$-balanced property.

Now we define the main quantity to be studied in the paper. A mean of $n$ variables on $I$ is called a \emph{generalized quasi-arithmetic mean} if 
there exist functions $f_1,\dots,f_n$ with the property that
\begin{itemize}
	\item[($\star$)] $f_1,\dots,f_n\colon I\to\R$ are continuous, monotone in the same 
	sense and not simultaneously constant on 
	any non-trivial interval in $I$,
\end{itemize}
such that the mean is of the form
\Eq{Mat}{
	M_f(x_1,\dots,x_n):=(f_1+\dots+f_n)^{-1}(f_1(x_1)+\dots+f_n(x_n))
}
for all $x_1,\dots,x_n\in I$, where $f:=(f_1,\dots,f_n)$. It is straightforward to check that this is a well-defined quantity. The 
coordinate-functions $f_1,\ldots,f_n$ are called the \emph{generating functions of $M_f$}. Obviously, $M_f$ reduces to a quasi-arithmetic mean if 
$f_1=\ldots=f_n$ holds. Results on means of the form \eq{Mat} were published and investigated first in \cite{BajPal09} and 
\cite{Mat10b}. Their equality problem was solved in \cite{MatPal15} and related theorems can be found in \cite{DarPal13}.

Below we formulate our main result.

\begin{theorem*}
Let $\sigma:\{1,\dots,n\}\to\{1,\dots,n\}$ be a bijection. Then the generalized quasi-arithmetic mean $M_f$ is 
$\sigma$-balanced if and only if it is a quasi-arithmetic mean.
\end{theorem*} 

\Rem{*}{
In view of this theorem and by the fact that the symmetric generalized quasi-arithmetic means are exactly the quasi-arithmetic ones, in the class of 
the former quantities, the $\sigma$-balancing property is equivalent to symmetry. As we now demonstrate, these two properties do not imply each other 
in general. On the one hand, it is easy to see that the coordinate means, that is, $M^1\colon (x,y)\mapsto x$ and $M^2\colon(x,y)\mapsto y\ ((x,y)\in 
I^2)$ are $\sigma$-balanced for any bijection $\sigma:\{1,2\}\to\{1,2\}$, but are not symmetic. On the other hand, for a fixed parameter $t\in[0,1]$, 
define the symmetric mean $M:I\times I\to\R$ by
\Eq{*}{
M(x,y):=t\min(x,y)+(1-t)\max(x,y).
}
A straightforward computation yields that this quantity is balanced if and only if $t\in\{0,\frac12,1\}$ holds, i.e., exactly when $M$ is the 
arithmetic mean or an extremal one.
}

We note that, in light of Theorem D of the paper \cite{MatPal15}, our main result has an equivalent reformulation in terms of functional 
equations. 
Before presenting it, we mention the following assertion.

\begin{lemma}[Theorem D, \cite{MatPal15}]\label{L:charquasi}
A generalized quasi-arithmetic mean $M_{f}$ is a quasi-arithmetic one if and only if there are constants $D_1,\ldots,D_n\in\R$ for which 
$f_k=f_1+D_k\ (k=1,\ldots,n)$.
\end{lemma}

The statement below shows that Theorem provides us with the general solution of a certain functional equation. It is an immediate consequence of the 
theorem and lemma above.

\Cor{*}{
Let $\sigma:\{1,\ldots,n\}\to \{1,\ldots,n\}$ be a bijection. Then the system $f=(f_1,\ldots,f_n)$ of functions with the 
property $(\star)$ fulfills the functional equation
\Eq{*}{M_f\big(M_f(u_{\sigma(1)}(x,M_f(x))),\dots,M_f(u_{\sigma(n)}(x,M_f(x)))\big)=M_f(x)\qquad (x\in I^n)}
exactly when there exist a continuous strictly monotone map $\varphi:I\to\R$ and constants $D_1,\ldots,D_n\in\R$ satisfying 
\Eq{*}{f_k=\varphi+D_k\qquad (k=1,\ldots,n).}
}

\section{Auxiliary assertions and the proof of the main result}

Keeping the notation and conditions introduced in the previous section, define
\Eq{*}{
F:=f_1+\ldots+f_n,
}
and, for any tuple $x=(x_1,\ldots,x_n)\in I^n$ and number $k=1,\ldots,n$, let
\Eq{E:vk}{
v_k(x):=M_f\big(u_k(x,M_f(x))\big)=M_f(x_1,\ldots,x_{k-1},M_f(x),x_{k+1},\ldots,x_n).
}
It is clear that without loss of generality, we may and do assume that the $f_k$-s are increasing, since otherwise we could consider the $-f_k$-s instead of 
them, which are clearly generating functions of the same generalized quasi-arithmetic mean $(k=1,\ldots,n)$. It is also easy to see that it follows 
from our conditions that $F$ is an invertible function.

To establish the main result of the paper, we will need a series of statements. The first one provides a necessary condition for the 
$\sigma$-balancing property of a generalized quasi-arithmetic mean in terms of its generating functions. 

\Prp{P:necessarycond}{
If $M_f$ is $\sigma$-balanced, then for all tuples $x\in I^n$, one has
\Eq{*}{
\sum\limits_{j=1}^n\sum\limits_{k=1,\, k\ne j}^nf_j(v_{\sigma(j)}(x))-f_k(v_{\sigma(j)}(x))=0.
}
}

\begin{proof}
Let $x=(x_1,\ldots,x_n)\in I^n$. For each $k=1,\ldots,n$, we compute
\Eq{E:Fvk}{
F(v_{\sigma(k)}(x))=f_{\sigma(k)}(M_f(x))+\sum\limits_{j=1,\, j\ne {\sigma(k)}}^nf_j(x_j).
}
Since $M_f$ is $\sigma$-balanced, $\sum\limits_{j=1}^nf_j(v_{\sigma(j)}(x))=\sum\limits_{j=1}^nf_j(x_j),$
so $\sum\limits_{j=1,\ j\ne \sigma(k)}^nf_j(x_j)=\sum\limits_{j=1}^nf_j(v_{\sigma(j)}(x))-f_{\sigma(k)}(x_{\sigma(k)})$. This yields that the sum of 
the right-hand sides of the equality \eq{E:Fvk} over $k=1,\ldots,n$ equals
\Eq{*}{
\sum\limits_{k=1}^n\Big(f_{\sigma(k)}(M_f(x))+\sum\limits_{j=1}^nf_j(v_{\sigma(j)}(x))-f_{\sigma(k)}(x_{\sigma(k)})\Big),
}
thus, by summing that equation for those $k$-s, we infer
\Eq{*}{
\sum\limits_{k=1}^nF(v_{\sigma(k)}(x))
&=\sum\limits_{k=1}^n\sum\limits_{j=1}^nf_j(v_{\sigma(j)}(x))+\sum\limits_{k=1}^nf_{\sigma(k)}(M_f(x))-\sum\limits_{k=1}^nf_{\sigma(k)}(x_{\sigma(k)})\\
&=n\sum\limits_{j=1}^nf_j(v_{\sigma(j)}(x))+F(M_f(x))-\sum\limits_{k=1}^nf_k(x_k).
}
Now observe that $F(M_f(x))-\sum\limits_{k=1}^nf_k(x_k)=0,$ hence it follows that 
$\sum\limits_{k=1}^nF(v_{\sigma(k)}(x))=n\sum\limits_{j=1}^nf_j(v_{\sigma(j)}(x)),$ implying that
\Eq{*}{
\sum\limits_{k=1}^n\sum\limits_{j=1}^nf_j(v_{\sigma(k)}(x))=\sum\limits_{j=1}^n\sum\limits_{k=1}^nf_j(v_{\sigma(j)}(x)).
} 
The latter equation immediately yields the statement of \prp{P:necessarycond}.
\end{proof}

In what follows, our goal is to show that as $x$ varies in $I^n$, the quantities $v_1(x),\ldots,v_n(x)$ may vary independently of each other. To this 
end, given arbitrary numbers $c_k\in I$, we are going to investigate the solvability of the system of equations
\Eq{E:mainsys}{
c_k=v_k(x)=M_f(x_1,\ldots,x_{k-1},M_f(x),x_{k+1},\ldots,x_n),\qquad(k=1,\ldots,n)	
}
in $x=(x_1,\ldots,x_n)\in I^n$. Observe that in terms of the new unknowns $y_k$ and objects $y,\ d_k$ and $z_k(x)$ defined by
\Eq{*}{
	y_k:=f_k(x_k),\qquad 
	y:=(y_1,\ldots,y_n),\qquad
	d_k:=F(c_k),\qquad
	z_k(y):=f_k(M_f(x))=f_k\Big(F^{-1}\Big(\sum\limits_{j=1}^ny_j\Big)\Big),
}
one has that $x$ is a solution of \eq{E:mainsys} exactly when $y$ fulfills
\Eq{E:redsys}{
	d_k=y_1+\ldots+y_{k-1}+z_k(y)+y_{k+1}+\ldots+y_n\qquad(y_k\in f_k(I),\, k=1,\ldots,n).
}

In the sequel, we use the notation
\Eq{*}{
g_k:=f_k\circ F^{-1}\qquad\text{and}\qquad \overline{x}:=\frac1n\sum_{k=1}^nx_k\qquad(x=(x_1,\ldots,x_n)\in\R^n).
}
The statement below gives a formula for the only possible solution of \eq{E:redsys}.

\Prp{P:solsys}{
For each $(d_1,\ldots,d_n)\in F(I)^n$, the unique solution of \eq{E:redsys} is
\Eq{E:solsys}{
	y_k=\frac1{n-1}\Big((n-2)(g_k(\overline{d})-d_k)+\sum\limits_{j=1,\ j\ne k}^nd_j-g_j(\overline{d})\Big),
}
provided that for all $k=1,\ldots,n$, this expression belongs to $f_k(I)$.
}

\begin{proof}
Let $y\in f_1(I)\times\dots\times f_n(I)$ be any vector and observe that
\Eq{E:gkzk}{
	\sum\limits_{k=1}^ng_k(u)=u,\quad
	z_j(y)=g_j\Big(\sum\limits_{k=1}^ny_k\Big),\quad 
	\sum\limits_{k=1}^nz_k(y)=\sum\limits_{k=1}^ny_k\quad
	\Big(u\in\sum\limits_{k=1}^nf_k(I),\ j=1,\ldots,n\Big).
}

Assume that $y$ satisfies \eq{E:redsys}. We are going to show that then \eq{E:solsys} holds. In order to do this, we rewrite the $k$th equation 
of \eq{E:redsys} in the form
\Eq{E:shortredsys}{
	z_k(y)+\sum\limits_{j=1,\ j\ne k}^ny_j=d_k,\qquad(k=1,\ldots,n).
}
By summing these relations for all those $k$-s and using \eq{E:gkzk}, we obtain $\sum\limits_{k=1}^ny_k+\sum\limits_{k=1}^n\sum\limits_{j=1,\ j\ne 
k}^ny_j=\sum\limits_{k=1}^nd_k$. We deduce that $\sum\limits_{j=1}^ny_j=\overline{d},$ which gives us that
\Eq{*}{
	z_k(y)=z_k=g_k(\overline{d}),
}
so it does not depend on $y$.

Now we infer that, by defining
\Eq{*}{
	d:=(d_1,\ldots,d_n)\qquad\text{and}\qquad z:=(z_1,\ldots,z_n),
}
the system \eq{E:shortredsys} is equivalent to $Ay=(d-z)^T$, where $A$ is the $n\times n$ matrix whose diagonal entries are $0$ and the others are 
$1$, further $^T$ denotes transpose. We are going to show that $A$ is invertible, and therefore the unique solution of the latter equation is 
$A^{-1}(d-z)^T$. In fact, $A=\mathds{1}-E$, where $\mathds{1}$ and $E$ stands for the $n\times n$ all-ones and identity matrix, respectively. Observe 
that $P:=\frac1n\mathds{1}$ is symmetric and 
idempotent, i.e., an orthogonal projection. It follows that the spectral decomposition of $A$, as a linear operator, is $A=(n-1)P-(E-P)$, so its spectrum is $\{-1,n-1\}$, and 
thus it is invertible. Moreover, we get that
$$A^{-1}=\frac1{n-1}P-(E-P)=\frac n{n-1}P-E=\frac1{n-1}\mathds{1}-E,$$
which is the product of $\frac1{n-1}$ and the matrix whose diagonal entries are $2-n$, and the others are 1. Finally, referring to the previous discussion, 
\eq{E:solsys} follows.

By what we have proved so far, we see that the only possible solution of \eq{E:redsys} is the one given by \eq{E:solsys}. To establish that if it belongs to $f_1(I)\times\dots\times f_n(I)$, then it indeed satisfies the mentioned system, suppose that $y$ is defined by \eq{E:solsys} and the $d_k$-s are chosen such that the right-hand side of that equation is in $f_k(I)\ (k=1,\ldots,n)$. Then the previous paragraph implies 
$\sum\limits_{j=1,\ j\ne k}^ny_j=d_k-g_k(\overline{d})$, therefore, by adding these relations for all such $k$-s, we obtain 
$(n-1)\sum\limits_{k=1}^ny_k=\sum\limits_{k=1}^nd_k-g_k(\overline{d}).$ It entails that
$z_i(y)=g_i\Big(\frac1{n-1}\Big(\sum\limits_{k=1}^nd_k-g_k(\overline{d})\Big)\Big)$, yielding, in virtue of \eq{E:gkzk}, that
\Eq{*}{
z_i(y)+\sum\limits_{j=1,\ j\ne 
i}^ny_j=d_i-g_i(\overline{d})+g_i\Big(\frac1{n-1}\Big(\sum\limits_{k=1}^nd_k-\sum\limits_{k=1}^ng_k(\overline{d})\Big)\Big)=d_i\qquad(i=1,
\ldots,n).
}
We conclude that \eq{E:redsys} holds for $y$ and then the proof is complete.
\end{proof}

Motivated by the statement above, for $k=1,\dots,n$, we define
\Eq{*}{
	\alpha_k:=\frac1{n-1}\Big((n-2)(g_k(\overline{d})-d_k)+\sum\limits_{j=1,\ j\ne k}^nd_j-g_j(\overline{d})\Big),
}
which is the $k$th coordiante of the only possible solution of \eq{E:redsys}.

In virtue of \prp{P:solsys}, it is clear that establishing conditions on the $d_k$-s that guarantee the inclusion $\alpha_k\in f_k(I)$ is crucial for 
our investigation. They can be obtained by giving appropriate bounds for $\alpha_k$. This is established in the next assertion.

\Lem{L:bound}{
	For any $(d_1,\ldots,d_n)\in F(I)^n$, let
	\Eq{*}{
		c_*:=\min\{F^{-1}(d_1),\ldots,F^{-1}(d_n)\}\qquad\text{and}\qquad 
		c^*:=\max\{F^{-1}(d_1),\ldots,F^{-1}(d_n)\}.
	}
	Then given an arbitrary $k=1,\ldots,n$, one has
	\Eq{*}{
		f_k(c_*)+(F(c_*)-F(c^*))\le\alpha_k\le 
		f_k(c^*)+(F(c^*)-F(c_*)).
	}
}

\begin{proof}
Let
\Eq{*}{
	c_k:=F^{-1}(d_k)\quad\text{and}\quad M:=F^{-1}(\overline{d})\qquad(k=1,\ldots,n).
}
Then observe that $M$ is a quasi-arithmetic mean of $c_1,\ldots,c_n$, therefore $c_*\le M\le c^*$. We compute
\Eq{*}{
	\alpha_k&=\frac1{n-1}\Big((n-2)(f_k(M)-F(c_k))+\sum\limits_{j=1,\ j\ne k}^nF(c_j)-f_j(M)\Big)\\
	&\ge\frac1{n-1}\Big((n-2)(f_k(c_*)-F(c^*))+\sum\limits_{j=1,\ j\ne k}^nF(c_*)-f_j(c^*)\Big)\\
	&=\frac1{n-1}\Big((n-2)f_k(c_*)-(n-2)F(c^*)+(n-1)F(c_*)-\sum\limits_{j=1}^nf_j(c^*)+f_k(c^*)\Big)\\
	&=\frac1{n-1}((n-1)f_k(c_*)+(n-1)F(c_*)-(n-1)F(c^*))=f_k(c^*)+\big(F(c_*)-F(c^*)\big),
}
which immediately entails the first inequality in the lemma. The second can be proved in the same way. 
\end{proof}

An immediate consequence of the next statement is that locally, the $v_k(x)$-s, defined in \eq{E:vk}, may vary independently of each other as $x$ 
varies in $I^n\ (k=1,\ldots,n)$.

\Prp{P:locsolv}{
	For any point $p\in I^{\circ}$, there is an open interval $U_p\subset I^{\circ}$, such that given arbitrary elements $c_1,\ldots,c_n\in U_p$, the 
	system \eq{E:mainsys} is solvable in $I^n$.
}

\begin{proof}
By \prp{P:solsys} and the observation preceding it, in order to prove the assertion, it is enough to show that for each element $p\in I^{\circ}$, one can 
find an open interval $U_p\subset I^{\circ}$ satisfying that given any tuple $(c_1,\ldots,c_n)\in U_p^n$ of points, with the notation $d_k:=F(c_k)$, the inclusion 
$\alpha_k\in f_k(I)$ holds for all $k=1,\ldots,n$.

To do this, let $k$ be such a number. Notice that the latter inclusion is valid if $f_k(I)=\R$. Otherwise pick real numbers $a<p<b$ for which 
$[a,b]\subset I^{\circ}$. Since $f_k$ is continuous, $f_k(I)$ is an interval, say $\langle a_k,b_k\rangle$. Define
\Eq{*}{
	\varepsilon_k:=\min\{f_k(a)-a_k,b_k-f_k(b)\}
}
with the convention that if $a_k=-\infty$ and $b_k=+\infty$, then the first and second element is missing from the latter set, respectively. Clearly, $\varepsilon_k>0$. By the continuity of $F$, it is uniformly continuous on $[a,b]$, hence there is a real number $r_k>0$ such that 
$|F(t)-F(s)|<\varepsilon_k$ for all $s,t\in [a,b]$ with $|t-s|<r_k$. Let
\Eq{*}{
	r:=\min\{r_1,\ldots,r_k\}.
}

Now choose an interval $p\in U_p=(\alpha,\beta)\subset[a,b]$ for which $\diam(U_p)<r$ and let $(c_1,\ldots,c_n)\in U_p^n$ be an $n$-tuple. Define
\Eq{*}{
c_*:=\min\{c_1,\ldots,c_n\}\qquad\text{and}\qquad
c^*:=\max\{c_1,\ldots,c_n\}.
}
In view of \lem{L:bound}, we have
\Eq{E:boundsalpha}{
	f_k(c_*)+(F(c_*)-F(c^*))\le\alpha_k\le f_k(c^*)+(F(c^*)-F(c_*)),
}
thus, by the conditions concerning the generating functions, we obtain that 
\Eq{E:1stupperboundalpha}{
	\alpha_k\le f_k(b)+(F(\beta)-F(\alpha)).
}
Since $\beta-\alpha=\diam(U_p)<r$, referring to the previous paragraph, we deduce that if $b_k<+\infty$, then
\Eq{E:uniform}{
	F(\beta)-F(\alpha)<\varepsilon_k\le b_k-f_k(b).
}
Inequalities \eq{E:1stupperboundalpha} and \eq{E:uniform} imply that in this case
\Eq{E:2ndupperboundalpha}{
	\alpha_k<b_k.
}
This holds in the case $b_k=+\infty$, too. If $a_k=-\infty$, then $\alpha_k>a_k$. Otherwise \eq{E:boundsalpha} yields 
that $\alpha_k\ge f_k(a)+(F(\alpha)-F(\beta))$. Referring to \eq{E:uniform}, we obtain that $F(\alpha)-F(\beta)>\max\{f_k(b)-b_k,a_k-f_k(a)\}$ 
with the convention that the first element of the set is missing in the case $b_k=+\infty$. The last two inequalities entail 
$\alpha_k>a_k$. Since we also have \eq{E:2ndupperboundalpha} and the equality $f_k(I)=\langle a_k,b_k\rangle$, we conclude that $\alpha_k\in 
f_k(I)$ for all $k=1,\ldots,n$. Finally, the assertion of \prp{P:locsolv} follows readily.
\end{proof}

It has the next immediate consequence.

\Cor{independencevk}{
	For any point $p\in I^{\circ}$, there is an open interval $U_p\subset I^{\circ}$, such that $U_p^n$ is contained in the range of 
	$(v_1,\ldots,v_n)\colon I^n\to\R^n$.
}

The next statement is of key importance in establishing that the differences of the generating functions of a $\sigma$-balanced generalized 
quasi-arithmetic mean are constants.

\Prp{L:locconstdiff}{
If $M_f$ is $\sigma$-balanced, then for any point $p\in I^{\circ}$, there are an open interval $U_p\subset I^{\circ}$ and constants 
$D_1,\ldots,D_n\in\R$ such that
\Eq{*}{f_k|_{U_p}=f_1|_{U_p}+D_k} holds for all $k=1,\ldots,n$.}

\begin{proof}
Let $p\in I^{\circ}$ be a point. Then by \prp{P:necessarycond} and the previous corollary, there is an open interval $p\in U_p\subset I^{\circ}$ for which
\Eq{*}{
	\sum\limits_{j=1}^n\sum\limits_{k=1,\ k\ne j}^n(f_j(t_j)-f_k(t_j))=0\qquad(t_j\in U_p,\, j=1,\ldots,n).
}
Fix $k\in\{1,\dots,n\}$ and $t_j$, where $j\in\{1,\ldots,n\}\setminus\{k\}$, and let $t_k$ run through $U_p$. Then, by the last equality, 
$(n-1)f_k(t_k)=C_k+\sum\limits_{j=1,\ j\ne k}^nf_j(t_k)$ holds with some number $C_k\in\R$. It yields that for all $x\in U_p$ and $k=1,\ldots,n$, one 
has $(n-1)f_k(x)=C_k+\sum\limits_{j=1,\ j\ne k}^nf_j(x)$. By subtracting the first equality from the $k$th in the last system of equations, we get
$$
(n-1)f_k(x)-(n-1)f_1(x)=C_k-C_1+f_1(x)-f_k(x)\quad(k=2,\ldots,n),
$$
which reduces to
\Eq{*}{
	f_k(x)=f_1(x)+\frac1n(C_k-C_1)\qquad(x\in U_p).
}
By setting $D_k:=\frac1n(C_k-C_1)$ for $k=1,\dots,n$, the assertion of \prp{L:locconstdiff} follows.
\end{proof}

Now we are in a position to verify the main result of the paper.

\begin{proof}[Proof of Theorem]
First, suppose that $M_f$ is a quasi-arithmetic mean. Then $f_1=\ldots=f_n=:\varphi$, and we have that $M_f$ is $\sigma$-balanced exactly when
$$
\frac1n\sum_{k=1}^n\varphi(x_k)=\frac1n\sum_{k=1}^n\varphi(v_k(x))\qquad(x=(x_1,\ldots,x_n)\in I^n).
$$
We compute
$$
\varphi(v_k(x))=\frac1n\Bigg(\varphi(M_f(x))+\sum_{j=1,\ j\ne k}^n\varphi(x_j)\Bigg)=\frac1n\Bigg(\frac1n\sum_{i=1}^n\varphi(x_i)+\sum_{j=1,\ j\ne 
k}^n\varphi(x_j)\Bigg),
$$
and now it follows that $M_f$ is $\sigma$-balanced if and only if
$$
\frac1n\sum_{k=1}^n\varphi(x_k)=\frac1n\sum_{k=1}^n\frac1n\Bigg(\frac1n\sum_{i=1}^n\varphi(x_i)+\sum_{j=1,\ j\ne 
k}^n\varphi(x_j)\Bigg)\qquad(x=(x_1,\ldots,x_n)\in I^n).
$$
This is \eq{Bal} in the setting where $\sigma$ is the identity, $M$ is the arithmetic mean and the $k$th argument is $\varphi(x_k)$. By the third 
paragraph of the introduction, the latter operation is $\sigma$-balanced, i.e., it satisfies \eq{Bal}. The previous discussion yields that $M_f$ is 
$\sigma$-balanced.

Now assume that $M_f$ has this property and let $p\in I^{\circ}$ be a given point. Then, by the last proposition, there exist an open interval $p\in 
U\subset I^{\circ}$ and numbers $D_k\in\R$ fulfilling the equalities
\Eq{E:repres}{
	f_k|_U=f_1|_U+D_k\qquad(k=1,\ldots,n).
}
It is obvious that we have a maximal open interval $p\in U_0\subset I^{\circ}$ for which the latter relations hold.
	
Now suppose for a moment that $U_0\ne I^{\circ}$. Then at least one endpoint of $U_0$ is in $I^{\circ}\cap\R$. Denote this element by $q$. Referring 
to \prp{L:locconstdiff}, there are an open interval $q\in V\subset I^{\circ}$ and numbers $D_k'\in\R$ satisfying $f_k|_V=f_1|_V+D_k'\ (k=1,\ldots,n)$. 
Since $U_0\cap V\ne\emptyset$ and \eq{E:repres} is valid for $U_0$, it holds with both constants $D_k,D_k'$ on that intersection, so $D_k=D_k'\ 
(k=1,\ldots,n)$. We conclude that $f_k$ can be represented in the form \eq{E:repres} on the open interval $U_0\cup V\subset I^{\circ}$ containing 
$U_0$ as a proper subset. However, this contradicts to the maximality of the latter interval. We conclude that $U_0=I^{\circ}$, which implies that 
$f_k$ has the form \eq{E:repres} on $I^{\circ}$. By the continuity of $f_k$, we infer that $f_k=f_1+D_k$ on $I$, and referring to 
Lemma \ref{L:charquasi}, this yields that $M_f$ is a quasi-arithmetic mean. The proof of our main Theorem is complete.
\end{proof}

\end{document}